\documentclass[11pt,a4paper,twoside]{article}
\usepackage{amsthm, amsfonts,amsmath}

\newtheorem{theorem}{Theorem}[section]
\newtheorem{lemma}{Lemma}[section]

\newtheorem{corollary}{Corollary}[section]

\theoremstyle{remark}
\newtheorem{remark}{Remark}[section]
\theoremstyle{remark}

\theoremstyle{definition}

\numberwithin{equation}{section}

\author{Vladimir Rovenski\footnote{University of Haifa, Dept. of Mathematics, 3498838, Haifa, Israel. e-mail: \texttt{vrovenski@univ.haifa.ac.il}}, 
Sergey Stepanov\footnote{Financial University, Dept. of Mathematics, 125993, Moscow, Russia; and 
\newline$\phantom.$\quad\
VINITI, Dept. of Mathematics, 20, 
125190 Moscow, Russia. e-mail: \texttt{s.e.stepanov@mail.ru}
}
and 
Irina Tsyganok\footnote{Financial University, Dept. of Mathematics, 125993, Moscow, Russia. e-mail: \texttt{i.i.tsyganok@mail.ru}}
}

\title{Back to almost Ricci solitons}

\begin{document}

\date{}

\maketitle

\begin{abstract}
	In the paper, we study complete almost Ricci solitons using the concepts and methods of geometric dynamics and geometric analysis.
In particular, we characterize Einstein manifolds in the class of complete almost Ricci solitons.
Then, we examine compact almost Ricci solitons using the orthogonal expansion of the Ricci tensor,
this allows us to substantiate the concept of almost Ricci solitons.

\vskip1.5mm
\noindent
\textbf{Keywords}: {Almost Ricci soliton; energy density; infinitesimal harmonic transformation; conformal Killing vector.}

\vskip.5mm
\noindent
\textbf{Mathematics Subject Classifications (2010)} Primary: 53C21; Secondary: 58J05.
\end{abstract}

%\begin{document}
%\maketitle

%% Citations in the text should be identified by appropriate numbers in square brackets, and consecutive references should be concatenated (e.g. [7, 12-15]).

%%% CONTENT

\section{Introduction}

One of the important components of the theory of Ricci flow are self-similar solutions called Ricci solitons, see \cite[pp.~153-176]{2}.
{Ricci solitons}, which are a generalization of Einstein manifolds, have been studied more and more intensively in the last twenty years.
This theory, besides being known after G.~Perelman's proof of the Poincar\'{e} conjecture (for details see \cite{1}),
has a wide range of applications in differential geometry and theoretical physics.
In turn, the study of almost Ricci solitons, which are a generalization of quasi-Einstein manifolds and Ricci solitons,
 was started by Pigola, Rigoli, Rimoldi, and Setty, see \cite{3}.
An $n$-dimensional $(n\ge 2)$ Riemannian manifold $(M,g)$ is called an \textit{almost Ricci soliton},
if there exist a smooth complete vector field $\xi$ and
a function
$\lambda\in C^\infty(M)$ such~that
\begin{equation}\label{GrindEQ__1_1_}
 {\rm Ric}=\frac{1}{2}\,{\cal L}_{\xi}\,g+\lambda\,g.
\end{equation}
Here, $Ric$ is the Ricci tensor and ${\cal L}_{\xi}$ is the Lie derivative operator in the direction of $\xi$.
Namely, $({\cal L}_{\xi}\,g)(X,Y) = g(\nabla_X \xi, Y) +g(\nabla_Y \xi, X)$ for all smooth vector fields $X,Y$ on~$M$,
where $\nabla$ is the covariant derivative (Levi-Civita connection).
Denote by $(M,g,\xi,\lambda)$ an {almost Ricci soliton}. For $\lambda=const$, it is a Ricci soliton.
Note that when $\xi$ is a Killing vector field, i.e., ${\cal L}_{\xi}\,g=0$, an almost Ricci soliton $(M,g,\xi,\lambda)$ is
\textit{Einstein manifold}, i.e., $Ric=\frac{s}{n}\,g$, from which we can apply Schur's lemma, e.g., \cite{6}, to obtain $\lambda=const$.
In the special case, where $\xi=\nabla f$ for some function $f\in C^\infty(M)$, we say that $(M,g,\xi,\lambda)$ is a gradient almost Ricci soliton with potential function $f$.

In \cite{3}, almost Ricci complete gradient solitons are considered.
Other more recent papers have studied compact almost Ricci solitons (e.g., \cite{4,19,20}) or almost Ricci solitons on manifolds with additional geometric structures, e.g., \cite{21,22}.
There are also attempts to find applications of almost Ricci solitons in theoretical physics, see, e.g., \cite{5}.

In Sections~\ref{sec:02}-\ref{sec:03}, we study complete almost Ricci solitons using concepts and methods of geometric dynamics and geometric analysis. In Section~\ref{sec:04}, we study compact almost Ricci solitons applying the orthogonal expansion
of symmetric two-tensors (see \cite[p.~130]{17}) to the Ricci tensor.
In~particular, this will make it possible to substantiate the concept of almost Ricci solitons.

\section{Complete almost Ricci solitons}
\label{sec:02}

Here, we study complete almost Ricci solitons from the point of view of geometric dynamics, see \cite{9,10}.
Denote by $\theta$ the $g$-dual one-form of $\xi$
and $\bar\Delta =\nabla^*\nabla$ the Laplace operator for the formal adjoint operator ${\nabla}^*$ of $\nabla$.
First, we formulate a lemma needed to prove our main results.

\begin{lemma}\label{L-01}
The vector field $\xi$ of an almost Ricci soliton $(M,g,\xi,\lambda)$ satisfies the equation
\begin{equation}\label{GrindEQ__1_3_}
 \bar\Delta\,\theta = {\rm Ric}(\xi,\cdot) - (n-2)\,d\lambda\,.
\end{equation}
\end{lemma}

Recall that a vector field $\xi$ generates a flow on a manifold, which is a one-parameter group of infinitesimal self-diffeomorphisms \cite[pp.~12-14]{6}.
A vector field $\xi$ is an \textit{infinite\-simal harmonic transformation} on $(M,g)$ if the local one-parameter group of infinitesimal self-diffeomorphisms ge\-nerated by $\xi$ is a group of harmonic self-diffeomorphisms.
A vector field $\xi$ is an infinitesimal harmonic transformation in $(M,g)$ if and only if
$\bar{\Delta}\,\theta=Ric(\xi,\,\cdot)$, where $\theta$ is the $g$-dual one-form of $\xi$, see \cite{8}.
In particular, the Killing vector field is an example of an infinitesimal harmonic transformation on $(M,g)$, see~\cite{7}.
Moreover, a~vector field $\xi$ associated with a Ricci soliton $(M,g,\xi,\lambda)$ is also an infinitesimal harmonic transformation on $(M,g)$ \cite{8,23}. Note that a local one-parameter group of infinitesimal harmonic transformations, or a harmonic flow generated by $\xi$, is directly related to De Turck harmonic flows \cite[pp.~113--117]{2}.

The next corollary follows from Lemma~\ref{L-01}.

\begin{corollary}\label{C-01} An $n$-dimensional $(n\ge 3)$ almost Ricci soliton $(M,g,\xi,\lambda)$ is a Ricci soliton if and only if its vector field $\xi$ is an infinitesimal harmonic transformation. At the same time, the vector field $\xi$ associated with a two-dimensional almost Ricci soliton $(M,g,\xi,\lambda)$ is an infinitesimal harmonic transformation.
\end{corollary}

The function
\begin{equation*}
 e(\xi):=\frac{1}{2}\,\|\xi\|^2=\frac{1}{2}\,g(\xi,\xi).
\end{equation*}
is said to be the \textit{energy density} of the flow generated by the vector field $\xi$, see \cite[pp.~273--274]{9}.
The \textit{kinetic energy} of the flow of $\xi$ is defined by the integral formula, see \cite[pp.~2; 19; 37]{10},
\begin{equation*}
 E(\xi)=\int_{M}\,e(\xi)\,d\,{\rm vol}_{g}\,.
\end{equation*}
The kinetic energy can be infinite or finite (e.g., on a compact manifold). Note that the kinetic energy plays an impotent role in Hamilton dynamics, see, e.g., \cite{10}.

  Based on the above definition and Lemma~\ref{L-01}, we formulate our main theorem.

\begin{theorem}\label{T-01}
Let $(M,g,\xi,\lambda)$ be an $n$-dimensional $(n \ge 3)$ complete almost Ricci soliton such that
the rate of change of $\lambda$ along the trajectories of the $\xi$-flow is pointwise bounded from below by ${\rm Ric}(\xi,\xi):$
\begin{equation}\label{Eq-A}
 {\cal L}_{\xi}\,\lambda \ge {\rm Ric}(\xi, \xi).
\end{equation}
If~$E(\xi)<\infty$, then $\xi$ is a parallel vector field and $(M,g)$ is an Einstein manifold.
Furthermore, if the soliton has infinite volume, then $\xi=0$.
\end{theorem}

According to the definition of the Lie derivative, e.g., \cite[pp.~29-30]{6}, the Lie derivative of a function $f\in C^1(M)$ with respect to vector field $\xi$ in \eqref{Eq-A} is given by 
\[
 {\cal L}_{\xi}f =\xi(f)=df(\xi)=\nabla_\xi f=g(\nabla f, \xi).
\]
In particular, it follows from Theorem~\ref{T-01} that not every complete Riemannian manifold supports an almost Ricci soliton structure, see also \cite[Corollary 1.5 and Example 2.4]{3}.

\begin{remark}\rm
%Recall that
If a Riemannian manifold $(M,g)$ admits a complete parallel vector field, then $(M,g)$ is reducible,
i.e., is locally the metric product of a real line and some other Riemannian manifold.
In~Theorem~\ref{T-01}, instead of condition ``infinite volume" one can assume that $(M,g)$ is not reducible.

Theorem~\ref{T-01} can be supplemented as follows.
Let $(M,g,\xi,\lambda)$ be a two-dimensional complete almost Ricci soliton satisfying $Ric(\xi, \xi)\le0$ and $E(\xi)<\infty$,
then $(M,g,\xi,\lambda)$ is isometric to Euclidean plane or one of the flat complete surfaces: cylinder, torus, M\"{o}bius band and Klein~bottle.
\end{remark}

 The following assertion follows from \eqref{GrindEQ__1_1_} and Theorem~\ref{T-01}.

\begin{corollary}\label{C-02}
Let $(M,g,\xi,\lambda)$ be a complete almost Ricci soliton such that $E(\xi)<\infty$.
If $\lambda$ is a non-decreasing function along trajectories of this flow and ${\cal L}_{\xi}\sqrt{e(\xi)}\le -\lambda$, then $(M,g)$ is an Einstein manifold. Furthermore, if the soliton has infinite volume, then $\xi=0$.
\end{corollary}

Recall that the \textit{volume form} of $(M,g)$ is defined by equality $\omega_g(\partial_1,\ldots,\partial_n)=\sqrt{\det g}$ for $\partial_k={\partial}/{\partial x^k}$ with respect to local coordinates $x^1,\dots, x^n$.
Note also that a Riemannian manifold $(M,g)$ has a (global) volume element if and only if $(M,g)$ is orientable, see \cite[p.~195]{11}.
A~volume form on a connected manifold $(M,g) $ has a single global invariant, namely the (overall) volume,
$\mathrm{Vol}(M,g)=\int_M {\omega_g}\,d\,{\rm vol}_{g}$, which is invariant under volume-form preserving transformations.
The volume $\mathrm{Vol}(M,g)$ can be infinite or finite (e.g., $\mathrm{Vol}(M,g)<\infty$ for a compact manifold~$M$).

On the other hand, a complete non-compact Riemannian manifold with non-negative Ricci curvature has infinite volume, see~\cite{12}.
For the volume form $\omega_g$ of $(M,g)$, one can consider its \textit{Lie derivative} along trajectories of the flow
of $\xi$. Namely, we have the following, see \cite[p.~281]{1}:
 ${\cal L}_\xi\,\omega_g=({\rm div}\,\xi)\,\omega_g$.
According to the definition of the Lie derivative, ${\cal L}_\xi\,\omega_g$ measures the rate of the change of the volume form $\omega_g$ under deformations determined by a one-parameter group of differentiable transformations (or a flow) generated by the vector field $\xi$.

%On the other hand, 
In the well-known monograph \cite[p.~195]{11} the function ${\rm div}\,\xi$ was called the \textit{logarithmic rate of change of volume} (or, in other words, \textit{rate of volume expansion}) under the flow generated by the vector field $\xi$. On the other hand, the condition ${\rm div}\,\xi=0$ is equal to ${\cal L}_{\xi}\,\omega_{g}=0$. This means that the one-parameter group of differentiable transformations leaves $\omega_{g}$ invariant or, in other words, the vector field $\xi$ is an infinitesimal automorphism of the volume form, see \cite[p.~6]{9}. In dynamic, such a vector field $\xi$ is said to be \textit{divergence-free} and the flow generated by it is said to be \textit{incompressible}, see \cite[p.~125]{11}. The geometric dynamics of divergence-free vector fields was studied in detail in the monograph~\cite{10}.

 The following proposition is true.

\begin{lemma}\label{L-02} Let $(M,g,\xi,\lambda)$ be a complete oriented almost Ricci soliton such that the length of $\xi$ is integrable.
If the logarithmic rate of volumetric expansion doesn't change sign on $M$ under deformations determined by the flow
of $\xi$, then \eqref{GrindEQ__1_1_} has the form
\[
 {\rm Ric} =\frac{1}{2}\,{\cal L}_{\xi}\,g+\ \frac{s}{n}\,g\,,
\]
where components of the right hand side are orthogonal to each other with respect to the standard pointwise scalar product.
On the other hand, if \eqref{GrindEQ__1_1_} is of the form indicated above, then the flow of $\xi$ is incompressible.
\end{lemma}

%The following corollary of Lemma~\ref{L-02} is true.

\begin{corollary}\label{C-03}
Let $(M,g,\xi,\lambda)$ be a complete oriented almost Ricci soliton such that the length of $\xi$ is integrable.
If the logarithmic rate of volumetric expansion doesn't change sign on $M$ under deformations determined by the flow of $\xi$, then $\xi$ satisfies the equation
\[
 \bar{\Delta}\,\theta =Ric(\xi,\cdot) -\frac{n-2}{n}\,ds.
\]
\end{corollary}

%2
\section{Proof of results in Section~\ref{sec:02}}
\label{sec:03}

\textbf{Proof of Lemma~\ref{L-01} and Corollary~\ref{C-01}}.
The equation \eqref{GrindEQ__1_1_} has the following form with respect to local coordinates $x^{1},\,\ldots,\,x^{n}$:
\begin{equation} \label{GrindEQ__2_1_}
R_{ij}=\frac{1}{2}\,{\cal L}_{\xi}\,g_{ij}+\lambda\,g_{ij},
\end{equation}
where $R_{ij}$, $g_{ij}$ and ${\xi}_i$ stand, respectively, for the components of the Ricci tensor $Ric$, the metric tensor $g$, and the components ${\xi}_i=g_{ij}\,{\xi}^j$ of $\xi$.
Also $\theta={\xi}^j$ is the one-form corresponding to $\xi$ under the duality defined by the metric $g$.
According the formula ${\cal L}_{\xi}\,g_{ij} =\nabla_i\,{\xi}_j+{\nabla}_j\,{\xi}_i$,
where ${\nabla}_i = {\nabla}_{\partial/\partial x^i}$, we obtain from \eqref{GrindEQ__2_1_} the equality
\begin{equation} \label{GrindEQ__2_2_}
 {\rm div}\,\xi:={\nabla}_i\,\xi^i = s-n\,\lambda
\end{equation}
for the scalar curvature $s$ of the metric $g$. Applying the operator ${\nabla}^i=g^{ij}{\nabla}_j$ to \eqref{GrindEQ__2_1_}, we find
\begin{equation} \label{GrindEQ__2_3_}
 \nabla^i\,\nabla_{i}\,{\xi}_j+{\nabla}^i{\nabla}_j\,{\xi}_i={\nabla}_j\,s -2\,{\nabla}_j\,\lambda\,.
\end{equation}
Using the contracted second Bianchi identity, see \cite{6}, and ${\nabla}_i\,\xi^i$ of \eqref{GrindEQ__2_2_}, we have
\[
 {\nabla}_i{\nabla}_j\,{\xi}^i={\nabla}_j{\nabla}_i\,{\xi}^i +R_{ij}\,{\xi}^i ={\nabla}_j\,s -n{\nabla}_j\,\lambda +R_{ij}\,{\xi}^i.
\]
Using this and noting that $\bar\Delta=-\nabla^i\,\nabla_{i}$
(the Laplace operator $\bar\Delta =\nabla^*\nabla$ and its expression in coordinates coincide, see \cite[Paragraph~1.55]{17}),
we rewrite \eqref{GrindEQ__2_3_} in the following form:
\begin{equation}\label{GrindEQ__2_4_}
 \bar\Delta\,{\xi}_j =R_{ij}\,{\xi}^i -(n-2)\,{\nabla}_j\,\lambda\,.
\end{equation}
In coordinate-free form, \eqref{GrindEQ__2_4_} coincides with \eqref{GrindEQ__1_3_} that proves Lemma~\ref{L-01}.
By this and \eqref{GrindEQ__1_3_}, we complete the proof of Corollary~\ref{C-01}.
\hfill$\square$

\smallskip

\textbf{Proof of Theorem~\ref{T-01}}.
In our case, the second \textit{Kato inequality}, $\|\xi\|\,\Delta\|\xi\| \ge -g(\bar{\Delta}\,\xi,\xi)$, see \cite[p.~380]{12},
\begin{equation*}
 \|\xi\|\,\Delta\|\xi\| \ge -g(\bar{\Delta}\,\xi,\xi),
\end{equation*}
using Lemma~\ref{L-01}, can be rewritten in the following form:
\begin{equation*}
 \|\xi\|\,\Delta\|\xi\| \ge (n-2)\,{\cal L}_{\xi}\lambda -Ric(\xi, \xi) ,
\end{equation*}
where $\Delta$ is the \textit{Laplace--Beltrami operator} defined by the equality $\Delta\,f={\rm trace}_g\nabla\,df$ for an arbitrary  function $f\in C^2(M)$. The assumption \eqref{Eq-A} for $n\ge3$ implies that $\|\xi\|\,\Delta\|\xi\|\ge0$.
Then by the classical theorem of geometric analysis (see \cite{13}), either
$\int_M \|\xi\|^p\,d\,{\rm vol}_{g}=\infty$ for a positive number $p>1$, or $\|\xi\|=const$.
Thus, if $\|\xi\|\in L^p(M,g)$ at least for one $p>1$, then $\|\xi\|=const$.
Note that the inequality $E(f)<\infty$ is equivalent to $\|\xi\|\in L^2(M,g)$.
By~the above (for $p=2$), and the condition $E(f)<\infty$, we get $\|\xi\|=const$.
Using this and \eqref{GrindEQ__2_4_}, we derive
\begin{equation}\label{GrindEQ__2_5_}
 0=\frac{1}{2}\,\Delta\,g(\xi,\xi)
 =-g(\bar{\Delta}\,\xi,\xi)+ \|\nabla\,\xi\|^2=-Ric(\xi,\,\xi) + (n-2)\,{\cal L}_{\xi}\,\lambda +\|\nabla\,\xi\|^2.
\end{equation}
By \eqref{GrindEQ__2_5_}, $\xi$ is a parallel vector field, in particular, ${\cal L}_{\xi}\,g = 0$.
From \eqref{GrindEQ__1_1_} we get $Ric=\lambda\,g$.
Since $n\ge3$, by Schur's lemma, e.g., \cite{6}, we get $\lambda=const$.
Thus, $(M,g)$ is an Einstein manifold.
Next, if $\mathrm{Vol}(M,g)=\infty$, then  using $E(f)<\infty$ and $\|\xi\|=const$, we get $\xi=0$.
\hfill$\square$

\textbf{Proof of Corollary~\ref{C-02}}.
From \eqref{GrindEQ__2_1_} we derive the following equation:
\begin{equation} \label{GrindEQ__2_7_}
 R_{ij}\,{\xi}^i{\xi}^j = {\xi}^i ({\nabla}_i\,{\xi}_j)\,{\xi}^j +\lambda\,\|\xi\|^2,
\end{equation}
where ${\xi}^i({\nabla}_i\,\xi_j){\xi}^j =\frac{1}{2}\,{\xi}^i{\nabla}_i({\xi}_j\,{\xi}^j)={\cal L}_{\xi}\,e(\xi)$.
Thus, we can rewrite \eqref{GrindEQ__2_7_} in the form
\begin{equation*}
 {\rm Ric}(\xi, \xi) =\lambda\,\|\xi\|^2 + {\cal L}_{\xi}\,e(\xi).
\end{equation*}
Therefore, the condition $Ric(\xi, \xi) \le 0$ is equivalent to the inequality
 ${\cal L}_{\xi}\sqrt{e(\xi)} \le -\lambda$.
From the above, the validity of Corollary~\ref{C-02} follows.
\hfill$\square$

\smallskip

\textbf{Proof of Lemma~\ref{L-02} and Corollary~\ref{C-03}}.
Recall the following theorem, see \cite{16}: Let $\xi$ be a smooth vector field on a complete oriented Riemannian manifold $(M,g)$
such that $\|\xi\| \in L^1(M,g)$ and ${\rm div}\,\xi$ does not change sign on $(M,g)$, then ${\rm div}\,\xi=0$ on $(M,g)$.
In~particular, if $\xi$ is the vector field of a complete, noncompact and oriented almost Ricci soliton $(M,g,\xi,\lambda)$, then from \eqref{GrindEQ__2_2_} we obtain $s=n\,\lambda$. In this case, \eqref{GrindEQ__1_1_} can be rewritten in the form
\begin{equation}\label{GrindEQ__2_9_}
 {\rm Ric} = \frac{1}{2}\,{\cal L}_{\xi}\,g + \frac{s}{n}\,g .
\end{equation}
Hence, $g(\frac{1}{2}\,{\cal L}_{\xi}\,g,\,\frac{s}{n}\,g) =\frac{s}{n}\,{\rm div}\,\xi=0$. Therefore, the terms of the right hand side of \eqref{GrindEQ__2_9_} are orthogonal to each other with respect to the pointwise scalar product.
 In turn, \eqref{GrindEQ__1_3_} takes the form $\bar{\Delta}\,\theta = {\rm Ric}(\xi,\,\cdot) + \frac{n-2}{n}\,ds$.
\hfill$\square$

\section{Compact almost Ricci solitons}
\label{sec:04}

Here, we study compact almost Ricci solitons using the orthogonal expansion of the Ricci tensor, obtained using the Becce expansion of the space of symmetric two-tensors, see \cite[p.~130]{16}.

Denote by $S^pM$ the space of symmetric covariant $p$-tensors on
a compact Riemannian manifold $(M,g)$,
and define the \textit{global scalar product} for any $\varphi,{\varphi}'\in S^pM$ by the formula
\begin{equation} \label{GrindEQ__3_1_}
 \langle \varphi,\ \varphi'\rangle =\int_M {g(\varphi,\,\varphi')}\,d\,{\rm vol}_{g}.
\end{equation}
Let ${\delta}^*:C^{\infty}(S^1M)\to C^{\infty}(S^2M)$ be the first-order differential operator defined by ${\delta}^*\theta =\frac{1}{2}\,{\cal L}_{\xi}\,g$ for any smooth one-form $\theta$ and its $g$-dual vector field $\xi$, see \cite[p.~117; 514]{17}.
Let also $\delta : C^{\infty}(S^2M)\to C^{\infty}(S^1M)$ be the formal adjoint operator for ${\delta}^*$,
which is called the {divergence of symmetric two-tensors}. In this case, $\langle \varphi,{\delta}^*\theta\rangle =\langle \delta \varphi,\theta\rangle$ is true for any $\varphi\in C^{\infty}(S^2M)$ and $\theta\in C^{\infty}(S^1M)$.
For a compact Riemannian manifold $(M,g)$, the algebraic sum $\mathrm{Im}\,{\delta}^* +C^{\infty}(M)\,\cdot g$ is closed in $S^2M$,
and the following decomposition is true:
\begin{equation}\label{GrindEQ__3_2_}
 S^2M=(\mathrm{Im}\,{\delta}^* +C^{\infty}(M)\,\cdot g)\oplus({\delta}^{-1}(0)\cap{\mathrm{trace}}^{-1}_g(0));
\end{equation}
furthermore, both factors in \eqref{GrindEQ__3_2_} are infinite-dimensional and orthogonal to each other with respect to the global scalar product \eqref{GrindEQ__3_1_}, see~\cite[p.~130]{17}.

\begin{lemma} Let $(M,g)$ be a compact $n$-dimensional Riemannian manifold and 
\[
 {\rm Ric}=\frac{1}{2}\,{\cal L}_{\xi}\,g +\lambda\,g +\varphi
\]
the orthogonal expansion $($with respect to the global scalar product$)$ of the Ricci tensor
for some vector field $\xi$ and divergence-free and trace-free symmetric two-form $\varphi$. Then
\begin{enumerate}
\item for $n\ge 3$, the vector field $\xi$ is an infinitesimal harmonic transformation if and only if $\lambda =const$,
and for $n=2$, the vector field $\xi$ is an infinitesimal harmonic transformation;

\item the assumptions $n\ge 3$ and $\int_M ({\cal L}_{\xi}\,s)\,d\,{\rm vol}_{g}\ge 0$ imply that $\varphi=Ric-\frac{s}{n}\,g$ for the scalar curvature $s=const$, and that $\xi$ is a conformal Killing vector field.
\end{enumerate}
\end{lemma}

\proof For the Ricci tensor, decomposition \eqref{GrindEQ__3_2_} has the form
\begin{equation} \label{GrindEQ__3_3_}
 {\rm Ric} = \big(\frac{1}{2}\,{\cal L}_{\xi}\,g +\lambda\,g\big) +\varphi
\end{equation}
for some divergence-free and trace-free tensor $\varphi\in C^{\infty}(S^2M)$ and function $\lambda\in C^{\infty}(M)$.
From \eqref{GrindEQ__3_3_} we get
\begin{equation*}
 -\delta\theta := {\rm div}\,\xi= s-n\,\lambda .
\end{equation*}
Applying the operator $\delta$ to \eqref{GrindEQ__3_3_}, we find, see also \eqref{GrindEQ__1_3_},
\begin{equation*}
 \bar{\Delta}\,\theta = {\rm Ric}(\xi,\,\cdot) - (n-2)\,d\lambda.
\end{equation*}
Therefore, for $n\ge 3$, $\xi$ is an infinitesimal harmonic transformation if and only if $\lambda =const$,
and for $n=2$, the vector field $\xi$ is an infinitesimal harmonic transformation.
%%%%%%%%%%%%%%%%%%%%%%%%
 On the other hand, using $n\,\lambda=\delta\,\theta+s$, we derive the following equalities:
\begin{eqnarray}
\nonumber
 0&&=\langle\varphi,\,{\delta}^*\theta +\lambda g\rangle = \langle {\rm Ric} -{\delta}^*\theta -\lambda\,g, {\delta}^*\theta +\lambda\,g\rangle \\
\nonumber
 &&=\langle {\rm Ric},{\delta}^*\theta\rangle -\langle{\delta}^*\theta, {\delta}^*\theta\rangle -\langle \lambda\,g,{\delta}^*\theta\rangle
 +\langle {\rm Ric},\,\lambda\,g\rangle -\langle{\delta}^*\theta,\lambda\,g\rangle -\langle\lambda\,g,\lambda\,g\rangle  \\
\nonumber
 &&= \langle\delta {\rm Ric},\theta\rangle -\langle{\delta}^*\theta,{\delta}^*\theta\rangle
 -2\langle\lambda\,g,{\delta}^*\theta\rangle +\int_M(\lambda\,s)\,d\,{\rm vol}_{g} -n\int_M {\lambda}^2\,d\,{\rm vol}_{g}  \\
\nonumber
 &&= -\frac{1}{2}\,\langle ds,\theta \rangle -\langle{\delta}^*\theta,{\delta}^*\theta\rangle
 +2\langle d\lambda,\theta\rangle +\int_M \lambda(s-n\lambda)\,d\,{\rm vol}_{g}  \\
\nonumber
 &&= -\frac{1}{2}\,\langle ds,\theta \rangle -\langle{\delta}^*\theta,{\delta}^*\theta \rangle
 +2\langle d\lambda,\theta\rangle -\int_M(\lambda\,\delta\,\theta)\,d\,{\rm vol}_{g} \\
 \label{GrindEQ__3_6_}
 &&
 = -\frac{1}{2}\,\langle ds,\theta\rangle -\langle {\delta}^*\theta,{\delta}^*\theta\rangle +\langle d\lambda,\theta\rangle .
\end{eqnarray}
 Hence,
\begin{equation}\label{GrindEQ__3_7_}
 n\,\langle d\lambda,\theta\rangle =\langle d(\delta\,\theta +s),\theta \rangle=\langle d\delta\theta,\theta\rangle
 +\langle ds,\theta \rangle=\langle \delta \theta,\delta \theta\rangle +\langle ds,\theta\rangle .
\end{equation}
Therefore, from \eqref{GrindEQ__3_6_} and \eqref{GrindEQ__3_7_} we derive
\begin{equation*}
 \frac{n-2}{2n}\int_M ({\cal L}_{\xi}\,s)\,d\,{\rm vol}_{g} =-\langle{\delta}^*\theta,{\delta}^*\theta \rangle
 +\frac{1}{n}\,\langle \delta \theta,\delta \theta\rangle \le 0,
\end{equation*}
because $\|\varphi\|^2\ge\frac{1}{n}\,{({\rm trace}_g\,\varphi)}^2$ for any covariant two-tensor $\varphi$.

The assumptions $n\ge 3$ and $\int_M ({\cal L}_{\xi}\,s)\,d\,{\rm vol}_{g}\ge0$, or, in particular, ${\cal L}_{\xi}\,s \ge0$, imply
\begin{equation} \label{GrindEQ__3_9_}
 \langle {\delta}^*\theta,{\delta}^*\theta \rangle -\frac{1}{n}\,\langle\delta \theta,\delta \theta\rangle =0.
\end{equation}
On the other hand, the following equality is valid:
\begin{equation}\label{GrindEQ__3_10_}
 \big\|\frac{1}{2}\,{\cal L}_{\xi}\,g -\frac{1}{n}\,({\rm div}\,\xi)\,g\big\|^2 = g({\delta}^*\theta,{\delta}^*\theta)
 -\frac{1}{n}\,(\delta\,\theta)^2.
\end{equation}
From \eqref{GrindEQ__3_9_} and \eqref{GrindEQ__3_10_} we find that
\[
 \frac{1}{2}\,{\cal L}_{\xi}\,g =\frac{1}{n}\,({\rm div}\,\xi)\,g,
\]
i.e., $\xi$ is a conformal Killing vector field.
In~this case, from \eqref{GrindEQ__3_3_} we deduce that 
\[
 \varphi ={\rm Ric} -\frac{s}{n}\,g,
\]
to which we can apply Schur's lemma and then conclude that $s=const$.
\hfill$\square$

\smallskip
A statement similar to the following corollary was proved in \cite{8} for Ricci solitons.

\begin{corollary} Let $(M,g,\xi,\lambda)$ be an $n$-dimensional $(n\ge3)$ compact almost Ricci soliton such that
\[
 \int_M ({\cal L}_{\xi}\,s)\,d\,{\rm vol}_{g}\ge0
\]
for its scalar curvature $s$. Then $(M,g)$ is isometric to a Euclidean $n$-sphere.
\end{corollary}

\proof
For $\varphi=0$, equations \eqref{GrindEQ__3_3_} have the form of almost Ricci soliton equations~\eqref{GrindEQ__1_1_}.
Thus, if $(M,g,\xi,\lambda)$ is a compact almost Ricci soliton such that
\[
 \int_M \,({\cal L}_{\xi}\,s)\,d\,{\rm vol}_{g} \ge 0,
\]
or, in particular, ${\cal L}_{\xi}\,s\ge 0$, 
then it is an Einstein manifold, and $\xi$ is a conformal Killing vector field.
By \cite[Corollary~4.5]{18}, the soliton is
%isometric to
a Euclidean
%$n$-dimensional
sphere.
\hfill$\square$

\begin{remark}\rm The concept of almost Ricci soliton was introduced in \cite{3} as a Riemannian manifold $(M,g)$ satisfying the equation 
\[
 {\rm Ric}+\frac{1}{2}\,{\cal L}_V\,g =\lambda\,g,
\]
where $\lambda\in C^\infty(M)$ and $V$ is a smooth vector field on $M$.
The above equation and \eqref{GrindEQ__1_1_} are equivalent. Namely, if we suppose that $\xi=-V$, then we derive \eqref{GrindEQ__1_1_} from the above equation. Thus, the inequalities ${\cal L}_{\xi}\,s \ge 0$ and $\int_M ({\cal L}_{\xi}\,s)\,d\,{\rm vol}_{g} \ge 0$
(with the scalar curvature $s$ of the metric $g$)
can be rewritten in the form
${\cal L}_{V}\,s\le0$ and $\int_M ({\cal L}_{V}\,s)\,d\,{\rm vol}_{g}\le 0$, respectively, see also \cite{4}, where inequalities are replaced by equalities. For example, if $s=const$ along the trajectories of the flow of
$\xi$ on a compact almost Ricci soliton $(M,g,\xi,\lambda)$ with a nonconstant function $\lambda$, then $(M,g)$ is isometric to a Euclidean $n$-sphere, see~\cite{4}.
\end{remark}

\end{document}